\documentclass[12pt,a4paper]{amsart}
\usepackage{latexsym,time}

\numberwithin{equation}{section}
\newtheorem{theorem}{Theorem}[section]

\newtheorem{proposition}[theorem]{Proposition}
\newtheorem{corollary}[theorem]{Corollary}

\theoremstyle{definition}

\theoremstyle{remark}

\newtheorem{example}[theorem]{Example}

\newcommand\Ext{\operatorname{Ext}}

\begin{document}
 
\author[P.~Schenzel]{Peter Schenzel}
\title[Descent from the form ring]
{Descent from the form ring and Buchsbaum rings}
\address{Martin-Luther-Universit\"at Halle-Wittenberg, 
Fachbereich Mathematik und Informatik, 
D --- 06 099 Halle (Saale), Germany}
\email{schenzel@mathematik.uni-halle.de}

\maketitle
\section{Introduction and Main Results}

One of the major problems in commutative algebra is to recover 
information about a commutative ring $A$ from known properties 
of the form ring $G := G_A(\mathfrak q) = \oplus_{n\geq 0} 
{\mathfrak q}^n/{\mathfrak q}^{n+1}$ with respect to some 
ideal $\mathfrak q$ of $A$. There are Krull's classical results  
saying that $A$ is an integral domain resp. a normal 
domain if $G$ is an integral domain resp. a normal domain. 
It follows from the work \cite{AA}, \cite{CN}, \cite{HR} that 
several other properties of a homological nature, like 
regularity, Cohen-Macaulayness, Gorensteinness etc., descend 
from $G$ to $A$. In this note we want to pursue this point of 
view further. To this end let $Q$ denote the homogeneous ideal 
of $G$ generated by all the inital forms of element of 
$\mathfrak q$. For our purposes here we investigate the local 
cohomology modules $H^{\bullet}_Q(G)$ and $H^{\bullet}_{\mathfrak q}(A)$ 
of $G$ with respect to $Q$ and of $A$ with respect to 
$\mathfrak q$ resp. For their definition and basic properties 
see \cite{aG67}. The first result concerns the descent of the 
finiteness from $H^i_Q(G)$ to $H^i_{\mathfrak q}(A)$.

\begin{theorem} \label{1.1}
Let $\mathfrak q$ be an ideal of a commutative Noetherian 
ring $A$. Assume $H^i_Q(G)$ is a finitely generated graded 
$G$-module for all $i<t$. Then $H^i_{\mathfrak q}(A), \enspace i<t$, 
is a finitely generated $A$-module.
\end{theorem} 

Furthermore, let $H^i(Q;G)$ and $H^i({\mathfrak q};A)$ denote the 
Koszul cohomology of $G$ with respect to $Q$ and of $A$ with respect 
to $\mathfrak q$ resp. Note that changing the basis yields isomorphic 
cohomology modules. Hence, it is not necessary to fix a basis in 
our notation. It is well-known, see e.g. \cite{aG67}, that there 
are canonical homomorphisms
$$
\begin{array}{lclcl}
f^i_G & : & H^i(Q;G) & \to & H^i_Q(G) \quad \text{and} \\\vspace*{.5pt}
f^i_A & : & H^i({\mathfrak q};A) & \to & H^i_{\mathfrak q}(A).
\end{array}
$$
Our next result concerns the descent of the surjectivity from 
$f^i_G$ to the surjectivity of $f^i_A$. Note that the surjectivity 
of $f^i_A$ is the crucial point in the investigation of local 
Buchsbaum rings, see \cite{jS80}. For a graded $G$-module $M$ 
let $[M]_k, \enspace k \in \mathbb Z$, denote its $k$-th 
graded piece.

\begin{theorem} \label{1.2}
Assume there are integers $k$ and $t$ such that
$$
\begin{array}{lcl}
[H^i_Q(G)]_{n-i} & = & 0 \enspace \text{for} \enspace n \not= k-1, 
k \enspace \text{and} \enspace 0 \leq i < t \enspace \text{ and } 
\\\vspace*{.2pt}   
[H^t_Q(G)]_{n-t} & = & 0 \enspace \text{for} \enspace n> k.
\end{array}
$$
Then $f^i_A$ is surjective for all $i<t$ if and only if $f^i_G$ is 
surjective for all $i<t$.
\end{theorem}

If $(A,{\mathfrak m})$ is a local Noetherian ring, then the 
Buchsbaum property does not descend from $G=G_A(\mathfrak m)$ to 
$A$. Counterexamples are given by Steurich (University of Essen), 
cf. \ref{4.3}. So our result yields as a corollary, cf. \ref{4.1}, 
that under certain additional assumptions on $G$ the basic ring 
$A$ is a Buchsbaum ring if $G$ is a Buchsbaum ring. 

The third result has supplementary character. It gives a 
sufficient condition for $f^i_G, \enspace i<t$, to be a 
surjective homomorphism. As shown in Example \ref{4.4} it is not 
necessary.

\begin{theorem} \label{1.3}
Let $t$ be an integer. Assume for each $0 \leq i<j<t$ and all 
integers $p$ and $q$ such that
\begin{displaymath}
[H^i_Q(G)]_{p-i} \not= 0 \enspace\text{and} \enspace 
[H^j_Q(G)]_{q-j} \not= 0 
\end{displaymath}
we have $p-q \not=1$. Then the canonical homomorphism
\begin{displaymath}
f^i_G : H^i(Q;G) \to H^i_Q(G)
\end{displaymath}
is surjective for $i<t$, provided $QH^i_Q(G)=0$ for $i<t$.
\end{theorem}

The proof of \ref{1.1} and \ref{1.2} is based on a spectral 
sequence technique of Serre \cite{jpS}, devised for passing 
from the form ring to the underlying ring, see 2. Another spectral 
sequence argument yields the proof of \ref{1.3}. The above 
results apply to Buchsbaum rings, for this see 4. We refer to \cite[5.6]{cW}
for an introduction and the baic results concerning spectral sequences.  

The author is grateful to R\"udiger Achilles for stimulating 
discussions during the preliminary 
draft of this note several years ago.

\section{Auxiliary Spectral Sequences}

A filtration of a ring $A$ is a decreasing sequence of ideals 
$({\mathfrak a}_n)_{n \in\mathbb Z}$ of $A$ such that 
${\mathfrak a}_m {\mathfrak a}_n \subseteq {\mathfrak a}_{m+n}$ a
nd ${\mathfrak a}_0=A$. A filtered $A$-module is an $A$-module 
$M$ with a decreasing sequence $(M_n)_{n \in\mathbb Z}$ 
of $A$-submodules of $M$ such that ${\mathfrak a}_m M_n 
\subseteq M_{m+n}$. Let $\mathfrak a$ be an ideal of $A$. 
A filtration $(M_n)_{n \in\mathbb Z}$ of $M$ is called essentially 
$\mathfrak a$-adic if ${\mathfrak a}M_n \subseteq M_{n+1}$ for 
all $n$ and ${\mathfrak a}M_n = M_{n+1}$ for all suffciently 
large $n.$, By virtue of \cite[(0.11.1.3)]{aG61} a filtration is 
called exhaustive (resp. co-discrete) provided $\cap _{n\in \mathbb Z} 
M_n = M$ (resp. there is an integer $m$ such that $M_m = M$).

\begin{proposition} \label{2.1}
Let $\mathfrak q$ be an ideal of a commutative Noetherian ring $A$. 
Then there exists a spectral sequence
\begin{displaymath}
E^{pq}_1 = [H^{p+q}_Q(G)]_p \Longrightarrow_p 
E^{p+q} = H^{p+q}_{\mathfrak q}(A)
\end{displaymath}
whose $E^{pq}_{\infty}$-term is the $p$-th component of the 
graded module associated to an essentially $\mathfrak q$-adic 
filtration of $H^{p+q}_{\mathfrak q}(A)$.
\end{proposition}

\begin{proof}
Based on Serre's technique \cite{jpS} Achilles and Avramov, 
see \cite{AA},  constructed a spectral sequence
\begin{displaymath}
E^{pq}_1 = [\Ext^{p+q}_G(G(M),G(N))]_p \Longrightarrow_p 
E^{p+q} = \Ext^{p+q}_A(M,N),
\end{displaymath}
where $M,N$ are finitely generated $A$-modules with $G(M), \enspace 
G(N)$ their form modules with respect to $\mathfrak q$. Set 
$N=A$ and $M=A/{\mathfrak q}^n$. Because $G(A/{\mathfrak q}^n) = 
G/Q^n, \enspace n \geq 0$, the $E_1$-term becomes 
\begin{displaymath}
E^{pq}_1 = [\Ext^{p+q}_G(G/Q^n,G)]_p.
\end{displaymath}
By virtue of the canonical homomorphism induced by
\begin{displaymath}
A/{\mathfrak q}^{n+1} \to A/{\mathfrak q}^n 
\enspace\text{and}\enspace G/Q^{n+1} \to G/Q^n
\end{displaymath}
both sides from a direct system. According to \cite[(0.11.1)]{aG61} 
we may form its direct limit, which yields the desired spectral sequence. 
Note that the filtration induced by the direct limit of the spectral 
sequences is in general not co-discrete. But it is always exhaustive. 
\end{proof}

The above proposition is the crucial point of our investigation. 
It is based on Serre's technique \cite{jpS}, devised for passing 
from the tangent cone of a variety to the variety itself. For 
further investigations see \cite{AA}. The next proposition concerns 
the Koszul cohomology. Its proof follows from \cite{jpS}, 
see also \cite{AA}.

\begin{proposition} \label{2.2}
Let ${\mathfrak q}, A,Q$, and $G$ as before. Then, there is a 
convergent spectral sequence 
\begin{displaymath}
E^{pq}_1 = [H^{p+q}(Q;G)]_p \Longrightarrow_p 
E^{p+q} = H^{p+q}({\mathfrak q};A)
\end{displaymath}
whose $E^{pq}_{\infty}$-term is the $p$-th component of the 
graded module associated to an essentially $\mathfrak q$-adic 
filtration of $H^{p+q}({\mathfrak q};A)$.
\end{proposition}

\noindent {\bf Proof of \ref{1.1}.}  
Because $H^i_Q(G), \enspace i<t$, is a finitely generated 
graded $G$-module it follows that $[H^i_Q(G)]_n =0$ for all 
$|n| \gg 0$ and $0 \leq i < t$. That is, the filtration of 
$H^i_{\mathfrak q} (A), \enspace i<t $, given in \ref{2.1}, 
is finite. Since all the $E^{pq}_{\infty}$-terms, $p+q<t$, are 
annihilated by $\mathfrak q$ we see that 
$H^i_{\mathfrak q}(A), \enspace i < t$, is annihilated by a 
power of $\mathfrak q$. According to \cite[Lemma 3]{gF}, 
this completes the proof. 
\hfill $\Box$ \vspace*{5pt}

A litte bit more is true, if we assume $H^i_Q(G), \enspace i<t$, 
a $G$-module of finite length.

\begin{corollary} \label{2.4}
Suppose $H^i_Q(G), \enspace i<t$, is a $G$-module of finite length. 
Then $H^i_{\mathfrak q}(A), \enspace i<t$, is an $A$-module of 
finite length and 
$$
L_G(H^i_Q(G)) \geq L_A(H^i_{\mathfrak q}(A)) \text{ for } i<t.
$$
\end{corollary} 

The proof of \ref{2.4} is similarly to the proof of \ref{1.1}. 
Hence we omit it. The particular case of \ref{2.4} for a local 
ring $A$ and an $\mathfrak m$-primary ideal has been shown in 
\cite[(4.2)]{pS84}, by a completely different technique.

As a further auxiliary result we shall use the following well-known 
spectral sequence, see e.g. \cite[Section 2]{pS95}. 

\begin{proposition} \label{2.5}
There is a convergent spectral sequence
\begin{displaymath}
E^{pq}_2 = H^p(Q;H^q_Q(G)) \Longrightarrow E^{p+q} = H^{p+q}(Q;G)
\end{displaymath}
for computing the Koszul cohomology.
\end{proposition}

\section{Proof of (1.2) and (1.3)}

\noindent{\bf Proof of \ref{1.2}.}
Firstly, we consider the spectral sequence given in \ref{2.1}. We 
claim that  $E^{pq}_1 = E^{pq}_{\infty}$ for $p+q<t$. In order to 
show this we consider the subsequent stages, i.e.
\begin{displaymath}
E^{p-r,q+r-1}_r \to E^{pq}_r \to E^{p+r,q-r+1}_r.
\end{displaymath}
Assume $E^{pq}_r \not= 0$ for some $p+q=i<t$. Then $E^{pq}_1 \not= 0$, 
i.e. $p=k-i-1$ or $p=k-i$ by virtue of the assumption. 
Now $E^{p+r,q-r+1}_r$ is a subquotient of 
$E^{p+r,q-r+1}_1 = [H^{i+1}_Q(G)]_{p+r} =0$ for $p=k-i-1$ or $p=k-i$. 
Therefore, if $E^{pq}_r \not= 0$ we have $E^{p+r,q-r+1}_r =0$. 
The same arguments yield $E^{p-r,q+r-1}_r =0$ if $E^{pq}_r \not= 0$. 
That is, $E^{pq}_1 = E^{pq}_{\infty}$ for all $p+q = i<t$. Thus, 
$H^i_{\mathfrak q}(A), \enspace i<t$, possesses a filtration of two 
terms. Hence, there is a short exact sequence
\begin{displaymath}
0 \to [H^i_Q(G)]_{k-i} \to H^i_{\mathfrak q}(A) \to 
[H^i_Q(G)]_{k-i-1} \to 0.
\end{displaymath}  
Secondly, we use the spectral sequence given in \ref{2.5} in order 
to show that
\begin{eqnarray*}
[H^i(Q;G)]_{n-i} & = & 0 \enspace\text{for} \enspace n \not= k-1, k 
\enspace\text{and} \enspace 0 \leq i<t \enspace\text{and} \\\vspace*{.2pt}
[H^t(Q;G)]_{n-t} & = & 0 \enspace\text{for} \enspace n>k.
\end{eqnarray*}
To this end we note that $E^{pq}_2 = H^p(Q;H^q_Q(G))$ is by definition 
a subquotient 
of $H^q_Q(G)(p)^{\binom{m}{p}}$, where $m$ denotes the number of 
generators of $Q$. Therefore, $E^{pq}_2$ vanishes for $p+q=i<t$ 
(resp. $p+q=t$) in all graded pieces $\not= k-i-1$, $k-i$ (resp. $>k-t$). 
Hence, the spectral sequence proves the claim. Thirdly, we apply 
this result in order to investigate the spectral sequence given in 
\ref{2.2}. Similarly as in the first part of the proof it yields a 
short exact sequence
\begin{displaymath}
0 \to [H^i(Q;G)]_{k-i} \to H^i({\mathfrak q};A) \to [H^i(Q;G)]_{k-i-1} \to 0
\end{displaymath}
for $i<t$. Therefore, the canonical homomorphisms $f^i_G$ and $f^i_A$ 
induce a commutative diagram with exact rows
$$
\begin{array}{ccccccccc}
0 & \to & [H^i(Q;G)]_{k-i} & \to & H^i({\mathfrak q};A) & \to & 
[H^i(Q;G)]_{k-i-1} & \to &  0 \\
& & \downarrow [f^i_G]_{k-i} & &\downarrow f^i_A & & 
\downarrow  [f^i_G]_{k-i-1} &  & \\
0 & \to & [H^i_Q(G)]_{k-i} & \to & H^i_{\mathfrak q}(A) & \to & 
[H^i_G(G)]_{k-i-1} & \to & 0
\end{array}
$$
for $i<t$, where $[f^i_G]_n$ denotes the $n$-th graded piece of 
$f^i_G.$ According to \ref{3.2} $[f^i_G]_{k-i}$ is always surjective. 
Hence, the snake lemma yields that $f^i_A$, $i<t$, is surjective 
if and only if $f^i_G$, $i<t$, is ssurjective. 
\hfill $\Box$

\begin{proposition} \label{3.2}
Suppose there exists an integer $k$ such that 
$$
[H^i_Q(G)]_{k+1-i} =0 \text{ and all } i \in \mathbb Z.
$$ 
Then the canonical homomorphism 
\begin{displaymath}
[f^i_G]_{k-i} : [H^i(Q;G)]_{k-i} \to [H^i_Q(G)]_{k-i}
\end{displaymath}
is surjective for all $i$.
\end{proposition}

\begin{proof} For this we have to investigate the spectral sequence 
given in \ref{2.5}. We claim 
$[E^{0i}_2]_{k-i} = [E^{0i}_{\infty}]_{k-i}$ for all $i$. 
To this end we are looking at the subsequent stages. Incoming 
$d$'s come from subquotients of $E^{-r,i+r-1}_2 = 0$. Outgoing 
$d$'s land in  subquotients of 
$$
[E^{r,i-r+1}_2]_{k-i} = [H^r(Q;H^{i-r+1}_Q (G))]_{k-i}
$$ 
which is a subquotient of $[H^{i-r+1}_Q (G)(r)^{\binom{m}{r}}]_{k-i} = 0$, 
by virtue of the assumption. Here $m$ denotes the number of generators 
of $Q$. Hence the claim is proved. Next we note that
\begin{displaymath}
[E^{0i}_2]_{k-i} = [H^0(Q;H^i_Q(G))]_{k-i} \simeq [H^i_Q(G)]_{k-i}
\end{displaymath}
because $[H^i_Q(G)]_{n-i} = 0$ for $n>k$ by the assumption. 
But now the module $[H^i(Q;G)]_{k-i}$ possesses a filtration 
whose associated $0$-th graded piece is $[E^{0i}_{\infty}]_{k-i},$ i.e.  
there is a canonical surjective mapping
\begin{displaymath}
[H^i(Q;G)]_{k-i} \to [H^i_Q(G)]_{k-i}
\end{displaymath}
for all $i$. By virtue of the functoriality of the spectral sequence it 
is noting else but $[f^i_G]_{k-i}$.
\end{proof} 

\noindent {\bf Proof of \ref{1.3}.}
We use the spectral sequence
\begin{displaymath}
E^{pq}_2 = H^p(Q;H^q_Q(G)) \Longrightarrow E^{p+q}(Q;G)
\end{displaymath}
given in \ref{2.5}. Because of $QH^i_Q(G) = 0, i<t$, it follows that 
$$ 
E^{pq}_2 = H^q_Q(G)(p)^{\binom{m}{p}} \text{ for } p+q<t. 
$$ 
Here $m$ denotes the number of generators of $Q$. In the subsequent 
stages we have
\begin{displaymath}
E^{p-r,q+r-1}_r \to E^{pq}_r \to E^{p+r,q-r+1}_r.
\end{displaymath}
Suppose that 
\begin{displaymath}
[E^{pq}_r]_n \not= 0 \text{ for some } p+q<t, 
\end{displaymath}
then $[E^{pq}_2]_n \not= 0$. Because $E^{p+r,q-r+1}_r$ resp. 
$E^{p-r,q+r-1}_r$ is derived from $H^{p+r}(Q;H^{q-r+1}_Q(G))$ 
resp. $H^{p-r}(Q;H^{q+r-1}_Q (G))$ it follows that  
\begin{displaymath}
[E^{p+r,q-r+1}_r]_n = 0 \text{ resp.} [E^{p-r,q+r-1}_r]_n =0 
\end{displaymath} 
by virtue of the assumption. That is, $E^{pq}_2 = E^{pq}_{\infty}$ 
for $p+q<t$. Because $H^i(Q;G)$ possesses a filtration whose 
associated $i$-th graded piece is $E^{0i}_{\infty}$ there is a 
canonical surjective homomorphism
\begin{displaymath}
H^i(Q;G) \to H^i_Q(G), \enspace i<t.
\end{displaymath}
By virtue of the functoriality of the considered spectral 
sequence it is nothing else but $f^i_G$, i.e. $f^i_G$, $i<t$, 
is surjective, as required. 
\hfill $\Box$

\section{Applications to Buchsbaum Rings}

In this section $(A,\mathfrak m)$ denotes a local Noetherian ring of 
dimension $d$ with $\mathfrak m$ its maximal ideal. Then $A$ is 
called a Buchsbaum ring if the difference 
$C(A) := L(A/\mathfrak q) - e({\mathfrak q};A)$ is an invariant 
of $A$ not depending on the choice of a parameter ideal 
$\mathfrak q$ of $A.$ 
Here $L(A/\mathfrak q)$ and $e({\mathfrak q};A)$ denote resp. 
the length of $A/\mathfrak q$ and the multiplicity of $A$ with 
respect to $\mathfrak q,$  see \cite{SV} for further details.
In his crucial paper \cite{jS80} St\"uckrad showed that 
$A$ is a Buchsbaum ring if and only if the canonical homomorphism 
$f^i_A : H^i({\mathfrak m};A) \to H^i_{\mathfrak m} (A)$, $i<d$, 
is surjective. Let $G=G_A(\mathfrak m)$ denote the associated 
form ring and $M=G_+$ the irrelevant maximal ideal of $G$. 
The graded ring $G$ is called a Buchsbaum ring if $G_M$ is 
a local Buchsbaum ring. Hence, our main results \ref{1.2} and 
\ref{1.3} apply to the situation of Buchsbaum rings.

\begin{corollary} \label{4.1}
Let $G=G_A(\mathfrak m)$ denote the associated graded form ring. 
Assume there is an integer $k$ such that 
\begin{eqnarray*}
[H^i_M(G)]_{n-i} = 0 & \text{for} & n \not= k-1, k 
\enspace\text{and}\enspace 0 \leq i < d \enspace\text{ and } \\\vspace*{.1pt}
[H^d_M(G)]_{n-d} = 0 & \text{ for } & n>k.
\end{eqnarray*}
Then $A$ is a Buchsbaum ring if and only if $G$ is a Buchsbaum ring. 
In this case $L(H^i_M(G)) = L(H^i_{\mathfrak m}(A))$ for all $0 \leq i<d$. 
\end{corollary}

\begin{proof}
Readily it follows from \ref{1.2} using the characterization of 
Buchsbaum rings in terms of the surjectivity of $f^i_A$ and $f^i_G$. 
The statement on the length of the local cohomology modules is clear 
by virtue of the short exact sequence given in the proof of \ref{1.2}. 
\end{proof}

The `only if' part of \ref{4.1} is one of the main results of 
Goto's paper \cite[Theorem (1.1)]{sG82}. His proof is completely different. 
It does not use Serre's spectral sequence technique for passing 
from the tangent cone to the ring.

\begin{corollary} \label{4.2}
Let $G$ be as in \ref{4.1}. Assume for each $0 \leq i<j<d$ and all 
integers $p$ and $q$ with
\begin{displaymath}
[H^i_M(G)]_{p-i} \not= 0 \enspace\text{and} \enspace 
[H^j_M(G)]_{q-j} \not= 0
\end{displaymath}
we have that $p-q \not= 1$. Then $G$ is a Buchsbaum ring, 
provided $MH^i_M(G) =0$ for $i<d$. 
\end{corollary}

The proof follows by virtue of \ref{1.3} accordingly to the 
fact that $G$ is a Buchsbaum ring if $f^i_G$, $i<d$, is 
surjective. Note that \ref{4.2} was shown independently by 
St\"uckrad \cite[Prop. 3.10]{SV} not using a spectral sequence 
technique. Particular cases of it were obtained in \cite{pS82}, resp. 
by Goto and others.

We conclude this section with two examples concerning the 
assumptions in \ref{1.2} and \ref{1.3}. The following example 
was given by Steurich, see also \cite[(4.10)]{sG82}.

\begin{example} \label{4.3}
The condition in \ref{1.2} is the best possible. Set 
\begin{displaymath}
A=k[|x,y,z|]/(x^2,xy,xz-y^r,y^{r+1},xz^2), \, r \geq 3 
\text{ an integer, }
\end{displaymath} 
where $k[|x,y,z|]$ denotes the formal power series ring over a  
field $k$. Note that $\dim A=1$. Then
\begin{displaymath}
G := G_A(\mathfrak m) = k[x,y,z]/(x^2,xy,,xz,y^{r+1},y^r z)
\end{displaymath} 
is a Buchsbaum ring, i.e. $f^0_G$ is surjective. Furthermore, 
\begin{displaymath}
\begin{array}{ccl}
H^0_M(G) & = & G/M(-1) \oplus G/M(-r), \\\vspace*{.2pt}
[H^1_M(G)]_n & = & 0 
\text{ for }  n>r-2 \text{ and } 
[H^1_M(G)]_{r-2} \not=  0.  
\end{array}
\end{displaymath}
On the other hand, $A$ is not a Buchsbaum ring, i.e. $f^0_A$ is 
not surjective.
\end{example}

While the finite length of $H^i_M(G), i \not= \dim G,$ is inherited 
to $H^i_{\mathfrak m}(A), i \not= \dim A,$ for the form ring 
$G = G_{\mathfrak m}(A)$ of a local ring $(A, \mathfrak m),$ see \ref{1.1}, 
this is not true in the Buchsbaum case. It does not hold even in 
the quasi-Buchsbaum case, where quasi-Buchsbaum means that for $i \not= 
\dim A$ the local cohomology is a vector space over $A/\mathfrak m.$ This 
follows because the local ring $(A, \mathfrak m)$ in \ref{4.3} is not 
quasi-Buchsbaum. So one might continue in order to improve the result 
in \ref{1.1} by taking into acount the more subtle behaviour of $k$-Buchsbaum 
rings. 

The next example shows that the assumption in \ref{1.3} is not 
necessary for $f_G$, $i<t$, surjective.

\begin{example} \label{4.4}
Let $R=k[x_1,\ldots,x_6]$ denote the polynomial ring over an 
infinite field $k$. Using a technique of Griffith and Evans 
in \cite{sG81}  Goto constructed examples of homogeneous prime 
ideals $P \subset R$ such that $R_1 := R/P$ is a 4-dimensional 
graded domain with
$$
\begin{array}{ccl}
H^i_M(R_1) & = & 0, \enspace i \not=1,4, \quad H^1_M(R_1) \simeq k(-2), 
\enspace\text{ and } \\\vspace*{.2pt} 
[H^4_M(R_1)]_n & = & 0 \enspace\text{ for all } n \geq 0.
\end{array}
$$
Let $R_2 =k[y_1,y_2,y_3]/F$, $F$ a homogeneous form of degree 3, 
be the coordinate ring of a plane cubic. Note that $R_2$ is a 
two-dimensional Cohen-Macaulay ring with $[H^2_M(R_2)]_n = 0$ 
for all $n \geq 1$ and $[H^2_M(R_2)]_0 = k$. Let $S$ denote 
the Segre product of $R_1$ and $R_2$ over $k$. Using the 
K\"unneth formula, as done in \cite[Section 5]{pS82}, we obtain
\begin{eqnarray*}
H^i_M(S) & = & 0, \enspace i \not= 1,2,5, \quad H^1_M(S) \simeq k^{10}(-2), 
\enspace\text{ and } \\
H^2_M(S) & \simeq & k(0).
\end{eqnarray*} 
Therefore, the assumptions of \ref{1.3} are not fulfilled because 
in this case $(1+2) - (2+0) = 1$. We show that $f^i_S$, $i<5$, 
is surjective. For $i \not= 0,3,4$ this is trivially true. 
Using \ref{2.5} it follows readily for $i=1$. Therefore, it is 
enough to show the surjectivity of $[f^2_S]_0$. The K\"unneth 
relations induce a commutative diagram
\begin{eqnarray*}
[R_1]_0 \otimes  [H^2(M;R_2)]_0 & \to & [H^2(M;S)]_0 \\
\downarrow [\text{id}]_0  \otimes [f^2_{R_2}]_0 & & \downarrow [f^2_S]_0 
\\\vspace*{.1pt}
[R_1]_0 \otimes [H^2_M (R_2)]_0 &\simeq & [H^2_M(S)]_0.
\end{eqnarray*}
Because $R_2$ is a Cohen-Macaulay ring with $[H^2_M(R_2)]_n = 0$ 
for all $n \geq 1$ it follows easily that 
$[f^2_{R_2}]_0 : [H^2(M;R_2)]_0 \to [H^2_M(R_2)]_0$ is an isomorphism. 
Hence, $[f^2_S]_0$ is surjective as required.
\end{example}

\end{document}